\documentclass[11pt,a4paper]{article}

\usepackage[a4paper,left=1.65cm,right=1.65cm,top=1.65cm,bottom=1.65cm]{geometry}
\usepackage{t1enc}
\usepackage[utf8]{inputenc}
\usepackage{amsthm,amssymb}

\usepackage[mathcal,mathscr]{eucal}
\usepackage{enumerate}
\usepackage[sort,nocompress]{cite}
\usepackage{indentfirst}
\usepackage{todonotes}
\usepackage{subcaption}
\usepackage{hyperref}
\hypersetup{
    colorlinks=true,       
    linkcolor=blue,        
    citecolor=red,         
    filecolor=magenta,     
    urlcolor=cyan,         
    linktocpage=true
}

\theoremstyle{plain}
\newtheorem{theorem}{Theorem}
\newtheorem{lemma}[theorem]{Lemma}
\newtheorem{corollary}[theorem]{Corollary}
\newtheorem{claim}[theorem]{Claim}

\theoremstyle{definition}

\usepackage{amsmath}

\begin{document}

\title{Supermodularity in Unweighted Graph Optimization II:  \break
Matroidal Term Rank Augmentation }

\author{Krist\'of B\'erczi\thanks{ MTA-ELTE Egerv\'ary Research Group,
Department of Operations Research, E\"otv\"os University, P\'azm\'any
P. s. 1/c, Budapest, Hungary, H-1117. e-mail:  {\tt berkri\char'100
cs.elte.hu .}} \ \ and  \ {Andr\'as Frank\thanks{MTA-ELTE Egerv\'ary
Research
Group, Department of Operations Research, E\"otv\"os University,
P\'azm\'any P. s. 1/c, Budapest, Hungary, H-1117. e-mail:  {\tt
frank\char'100 cs.elte.hu .} } } }

\maketitle

\begin{abstract} Ryser's max term rank formula with graph theoretic
terminology is equivalent to a characterization of degree sequences of
simple bipartite graphs with matching number at least $\ell$.  In a
previous paper \cite{Berczi-Frank16a} by the authors, a generalization
was developed for the case when the degrees are constrained by upper
and lower bounds.  Here two other extensions of Ryser's theorem are
discussed.  The first one is a matroidal model, while the second one
settles the augmentation version.  In fact, the two directions shall
be integrated into one single framework.  \end{abstract}

\section{Introduction}

Ryser \cite{Ryser58} derived a formula for the maximum term rank of a
$(0,1)$-matrix with specified row- and column-sums.  In graph
theoretic terms, his theorem is equivalent to a characterization for
the existence of a degree-specified simple bipartite graph (bigraph
for short) with matching number at least $\ell.$ Several natural
extensions, like the min-cost and the subgraph version, turned out to
be {\bf NP}-hard, but in a previous paper \cite{Berczi-Frank16a},
we could extend Ryser's theorem to the degree-constrained case when,
instead of exact degree-specifications, lower and upper bounds are
imposed on the degrees of the bigraph.  An even more general problem
was also solved when, in addition, lower and upper bounds were imposed
on the number of edges.  The main tool in \cite{Berczi-Frank16a} for
proving these extensions was a general framework for covering an
intersecting supermodular function by degree-constrained simple
bipartite graphs.

In the present paper we consider two other extensions of Ryser's
theorem:  the augmentation and the matroidal version.  In the first
one, a given initial bigraph is to be augmented to get a simple
degree-specified bigraph with matching number at least $\ell.$ In
original matrix terms, this means that some of the entries of the
$(0,1)$-matrix are specified to be 1. The solvability of this version
is in sharp contrast with the {\bf NP}-completeness of another
variation when some entries of the matrix are specified to be 0. (This
follows from the {\bf NP}-completeness of the problem that seeks to
decide whether an initial bigraph $G_0$ has a perfectly matchable
degree-specified subgraph, see \cite{KaMk}, \cite{Palvolgyi},
\cite{Puleo}.)

In the matroidal extension of Ryser's theorem, there is a matroid on S and there is a matroid on T, and the goal is to find a degree-specified simple
bigraph including a matching that covers bases in both matroids.
These results will be consequences of a general framework including
both the augmentation and the matroidal cases.

The starting point in deriving the main result is the supermodular
arc-covering theorem by Frank and Jord\'an \cite{FrankJ31} (Theorem
\ref{thm:Frank-Jordan} below).  Since \cite{FrankJ31} describes a
polynomial algorithm, relying on the ellipsoid method, to compute the
optima in question, our matroidal term rank augmentation problem also
admits a polynomial algorithm.  One of the most important applications
in \cite{FrankJ31} is the directed node-connectivity augmentation
problem.  V\'egh and Bencz\'ur \cite{Benczur-Vegh} developed for this
special case a pretty intricate but purely combinatorial algorithm
(not relying on the ellipsoid method).  Although not mentioned
explicitly in \cite{Benczur-Vegh}, their algorithm can probably be extended to
work on the supermodular arc covering theorem when the function in question is $ST$-crossing supermodular, but the details have not been worked out.
(In the special case of node-connectivity augmentation, this general
oracle was realized via network flow computations.) Therefore the algorithm of V\'egh and Bencz\'ur seems to be adaptable to the term rank problem, too. In a forthcoming paper \cite{Berczi-Frank17}, we shall develop a much simpler algorithm along with a natural unification of the matroidal augmentation and the degree-constrained term rank problems.

\subsection{Notions and notation}

We use the notation of \cite{Berczi-Frank16a}.  Here we briefly repeat
the most important notions.  For a family $\cal T$ of sets, let $\cup
{\cal T}$ denote the union of the members of $\cal T$.  For a
subpartition ${\cal T}=\{T_1,\dots ,T_q\}$, we always assume that its
members $T_i$ are non-empty but $\cal T$ is allowed to be empty (that
is, $q=0$).

An arc $st$ {\bf enters} or {\bf covers} a set $X$ if $s\not \in X$,
$t\in X$.  A digraph {\bf covers} $X$ if it contains an arc covering
$X$.  Let $S$ and $T$ be two non-empty subsets of a ground-set $V$.
By an {\bf $ST$-arc} we mean an arc $st$ with $s\in S$ and $t\in T$.   Two sets
$X$ and $Y$ are {\bf $ST$-independent} if $X\cap Y\cap T=\emptyset $
or $S-(X\cup Y)=\emptyset $, that is, no $ST$-arc enters both sets.
Two subsets $X$ and $Y$ are {\bf comparable} if $X\subseteq Y$ or
$Y\subseteq X$.  Two non-comparable sets $X$ and $Y$ are {\bf
$T$-intersecting} if $X\cap Y\cap T\not =\emptyset $ and {\bf
$ST$-crossing} if $X\cap Y\cap T\not =\emptyset $ and $S-(X\cup Y)\not
=\emptyset $. A set-function $p$ is called {\bf positively $T$-intersecting ({\bf
$ST$-crossing}) supermodular} if the supermodular inequality
$$p(X)+p(Y)\leq p(X\cap Y)+p(X\cup Y)$$ holds for $T$-intersecting (resp. $ST$-crossing) subsets $X$ and $Y$ for which $p(X) > 0$ and $p(Y) > 0$. The function is \textbf{fully supermodular} if the supermodular inequality holds for every pair $X$ and $Y$ of subsets.

For a function $m:V\rightarrow {\bf R},$ the set-function $\widetilde
m$ is defined by $\widetilde m(X)=\sum [m(v):v\in X]$.  A set-function
$p$ can analogously be extended to families $\cal F$ of sets by
$\widetilde p({\cal F)}= \sum [p(X):X\in {\cal F}]$.

The following min-max theorem of Frank and Jord\'an \cite{FrankJ31}
will be a basic tool in the proof of the main theorem.

\begin{theorem}[Supermodular arc-covering, set-function version]
\label{thm:Frank-Jordan} A positively $ST$\--cross\-ing supermodular
set-function $p$ for which $p(V')\leq 0$ holds when no $ST$-arc enters
$V'$ can be covered by $\gamma $ (possibly parallel) $ST$-arcs if and
only if $\widetilde p({\cal I})\leq \gamma $ holds for every
$ST$-independent family $\cal I$ of subsets of $V$.  There is an
algorithm, which is polynomial in $\vert S\vert +\vert T\vert $ and in
the maximum value of $p(X)$, to compute the minimum number of
$ST$-arcs to cover $p$ and an $ST$-independent family $\cal I$ of
subsets maximizing $\widetilde p({\cal I})$.  \end{theorem}

Henceforth we assume that $S$ and $T$ are two disjoint non-empty sets
and $V:=S\cup T$.  Let $G\sp *=(S,T;E\sp *)$ denote the complete
bipartite graph on bipartition $(S,T)$.  Let $D\sp *=(S,T;A\sp *)$ be
the digraph arising from $G\sp *$ by orienting each of its edges from
$S$ to $T$, that is, $A\sp *$ consists of all $ST$-arcs.  More
generally, for a bigraph $H=(S,T;F)$, let $\overrightarrow
{H}=(S,T;\overrightarrow {F})$ denote the digraph arising from $H$ by
orienting each of its edges from $S$ toward $T$.

Throughout we are given a simple bigraph $H_0=(S,T;F_0)$ serving as an
initial bigraph to be augmented.  For $E_0:=E\sp *-F_0$, the bigraph
$G_0=(S,T;E_0)$ is called the {\bf bipartite complement} of $H_0$,
that is, $F_0$ and $E_0$ partition $E\sp *$.  Note that a bigraph
$G=(S,T;E)$ is a subgraph of $G_0$ precisely if the augmented bigraph
$G\sp +=(S,T;F_0+E)$ is simple.  For $X\subseteq S$ and $Y\subseteq
T$, let $d_{G_0}(X,Y)$ denote the number of edges of $G_0$ connecting
$X$ and $Y$.

\medskip 

\section{Matroidal covering and augmentation }

Let $p_T$ be a positively intersecting supermodular set-function on
$T$.  In \cite{Berczi-Frank16a}, we studied the problem of finding a
simple degree-specified bigraph $G=(S,T;E)$ covering $p_T$ in the
sense that $$\vert \Gamma _G(Y)\vert \geq p_T(Y) \hbox{ for every
subset}\ Y\subseteq T$$ where $\Gamma _G(Y)$ denotes the set of
neighbours of $Y$.  Here we consider a framework which is more general
in two directions.  First, for a given initial simple bigraph
$H_0=(S,T;F_0)$, we want to find a degree-specified bigraph $G$ in
such a way that $G\sp +:=G+H_0$ is simple and covers $p_T$.  This kind
of problems is often referred to as augmentation problems to be
distinguished from the synthesis problems where $F_0$ is empty.  If
$p_T\equiv 0,$ the augmentation problem is equivalent to finding a
degree-specified subgraph of the bipartite complement of $H_0$.

Second, we extend the notion of covering to matroidal covering in the
following sense.  Let $M_S=(S,r_S)$ be a matroid on $S$ with rank
function $r_S$.  A bigraph $G$ is said to {\bf $M_S$-cover} $p_T$ if
\begin{equation}r_S(\Gamma _G(Y)) \geq p_T(Y) \hbox{ for every subset $Y\subseteq
T$.}\ \label{eq:matfed} \end{equation} Clearly, when $M_S$ is the free matroid,
we are back at the original notion of covering by a bigraph.

\subsection{Degree-specified matroidal augmentation}

Let $m_V=(m_S,m_T)$ be a degree-spe\-ci\-fication.  A bigraph $G=(S,T;E)$
is said to {\bf fit} $m_V$ if $d_G(v)=m_V(v)$ for every $v\in S\cup
T$.  Our main goal is to describe a characterization for the existence
of a bigraph $G$ fitting $m_V$ so that $G+H_0$ is simple and
$M_S$-covers $p_T$.  The more general problem, when there are upper
and lower bounds on $V$, will be discussed in \cite{Berczi-Frank17}.
This degree-constrained version was solved in
\cite{Berczi-Frank16a} in the special case when $H_0$ has no edges and
$M_S$ is the $\ell$-uniform matroid on $S$.

Our main result is as follows.

\begin{theorem}\label{thm:msmt.nov} We are given a simple bigraph
$H_0=(S,T;F_0)$, a matroid $M_S=(S,r_S)$, a positively intersecting
supermodular set-function $p_T$ on $T$, and a degree-specification
$m_V=(m_S,m_T)$ on $V:=S \cup T$ for which $\widetilde m_S(S)=\widetilde m_T(T)=\gamma
$. There is a bigraph $G=(S,T;E)$ fitting $m_V$ for which $G\sp
+=G+H_0$ is simple and $M_S$-covers $p_T$ if and only if \begin{eqnarray}&\widetilde m_S(X)+ \widetilde m_T(Y) - d_{G_0}(X,Y)
+ \sum _{i=1}\sp q [p_T(T_i)- r_S(X\cup \Gamma _{H_0}(T_i))] \leq
\gamma & \nonumber \\& \hbox{whenever $Y\subseteq T$, $X\subseteq S$, and ${\cal
T}= \{T_1,\dots ,T_q\}$ is a subpartition of $T-Y,$ }\ &  \label{eq:ftgs} \end{eqnarray} where $G_0$ is the bipartite complement
of $H_0$.  \end{theorem}

\proof{Proof.} Necessity.  Suppose that there is a requested bigraph $G=(S,T;E)$
and let $G\sp +=(S,T;E \cup F_0)$.  Note that the simplicity of $G\sp
+$ is equivalent to the requirement that $G$ is a subgraph of $G_0$.
Let $X\subseteq S$ and $Y\subseteq T$ be subsets and let $\{T_1,\dots
,T_q\}$ be a subpartition of $T-Y$.  Let $W_i:=\Gamma _G(T_i) - [X\cup
\Gamma _{H_0}(T_i)]= \Gamma_{G^+}(T_i) - [X \cup \Gamma_{H_0}(T_i)].$ Then we have $$p_T(T_i) \leq r_S(\Gamma _{G\sp
+}(T_i)) \leq r_S(\Gamma _{H_0}(T_i)\cup X) + r_S(W_i) \leq $$ $$
r_S(\Gamma _{H_0}(T_i)\cup X) + \vert W_i\vert \leq r_S(\Gamma
_{H_0}(T_i)\cup X) + d_G(T_i,W_i)$$ from which $d_G(T_i,W_i)\geq
p_T(T_i) - r_S(X\cup \Gamma _{H_0}(T_i))$.  Therefore $G$ has at least
$\sum _{i=1}\sp q [ p_T(T_i) - r_S(X\cup \Gamma _{H_0}(T_i)) ]$ edges
connecting $T-Y$ and $S-X$, and $G$ has at least $ \widetilde m_S(X) +
\widetilde m_T(Y) - d_{G_0}(X,Y)$ edges with end-nodes in $X$ or in
$Y$, from which the inequality in \eqref{eq:ftgs} follows.

Sufficiency.  \ Let ${\cal H}_0:=\{V'\subseteq V:$ no arc of
$\overrightarrow {H_0}$ enters $V'\}$.  Then ${\cal H}_0$ is closed
under taking union and intersection.  In the following definition of set-function $p_0$, we have $X\subseteq S$, $Y\subseteq T$, and $y\in T$.

\begin{equation}p_0(V'):= \begin{cases}
\max\{p_T(y) - r_S(X),\ m_T(y) - \vert X\vert +
d_{H_0}(y) \} & \text{if $V'=X+y\in {\cal H}_0$,}\\ p_T(Y)-r_S(X) & \text{if $V'=X\cup Y\in {\cal H}_0$, $|Y|\geq 2$}\\
0 & \text{otherwise.} \end{cases}
\label{eq:p0def} \end{equation}

The definition of $p_0$ implies that $p_0(V')$ can be positive only if $V'\in\mathcal{H}_0$.

\begin{lemma}\label{lem:ujterm.Tsuper} The set-function $p_0$ is positively
$T$-intersecting supermodular.  \end{lemma}

\proof{Proof.} Let $X_1,X_2$ be subsets of $S$ and let $Y_1, Y_2$ be subsets
of $T$ for which $Y_1\cap Y_2\not =\emptyset $. Suppose that
$p_0(V_i)>0$ for $V_i=X_i \cup Y_i$ $(i=1,2)$.  Then each of the sets
$V_1$, $V_2,$ $V_1\cap V_2,$ and $V_1\cup V_2$ belongs to ${\cal
H}_0$.  We distinguish three cases.\medskip

\noindent {\bf Case 1} \ $p_0(V_i) = p_T(Y_i) -r_S(X_i)$ for $i=1,2$.
Then $$p_0(V_1) + p_0(V_2) =
[p_T(Y_1) -r_S(X_1)] + [p_T(Y_2) -r_S(X_2)] \leq $$ $$p_T(Y_1\cap Y_2)
-r_S(X_1\cap X_2) + p_T(Y_1\cup Y_2) -r_S(X_1\cup X_2) \leq p_0(V_1\cap
V_2) + p_0(V_1\cup V_2).$$

\noindent {\bf Case 2} \ $p_0(V_i) > p_T(Y_i) -r_S(X_i)$ for $i=1,2$.
Then $Y_1=Y_2 = \{y\}$ for some $y\in T$, and $p_0(V_i)= m_T(y) -
\vert X_i\vert + d_{H_0}(y)$.  We have $$p_0(V_1) + p_0(V_2)= m_T(y) -
\vert X_1\vert + d_{H_0}(y) + m_T(y) - \vert X_2\vert + d_{H_0}(y) =$$
$$m_T(y) -\vert X_1\cap X_2\vert + d_{H_0}(y) + m_T(y) - \vert X_1\cup
X_2\vert + d_{H_0}(y) \leq p_0(V_1\cap V_2) + p_0(V_1\cup V_2).$$

\noindent {\bf Case 3} \ $p_0(V_1) = p_T(Y_1) -r_S(X_1)$ and $p_0(V_2)
> p_T(Y_2) -r_S(X_2)$.  (The situation is analogous when the indices
$i=1,2$ are interchanged.)  Then $Y_2=\{y\}$ for some $y\in T$ and
$y\in Y_1$.  Since $$r_S(X_1\cup X_2) - r_S(X_1) \leq \vert (X_1\cup
X_2)-X_1\vert = \vert X_2\vert -\vert X_1\cap X_2\vert ,$$ we have
$-r_S(X_1) - \vert X_2\vert \leq -r_S(X_1\cup X_2) - \vert X_1\cap
X_2\vert $ and hence $$p_0(V_1) + p_0(V_2)= [p_T(Y_1)- r_S(X_1)] +
[m_T(y) - \vert X_2\vert + d_{H_0}(y) ] = $$ $$[p_T(Y_1\cup Y_2) -
r_S(X_1)] + [m_T(y) - \vert X_2\vert + d_{H_0}(y) ] \leq $$ $$
p_T(Y_1\cup Y_2) - r_S(X_1\cup X_2) + m_T(y) - \vert X_1\cap X_2\vert
+ d_{H_0}(y) \leq p_0(V_1\cup V_2) + p_0(V_1\cap V_2),$$ as required.
$\bullet$\endproof

\begin{claim}\label{claim:gsmin} $m_S(s) \leq d_{G_0}(s)$ for each $s\in S$.  \end{claim}

\proof{Proof.} By applying \eqref{eq:ftgs} to $Y=T$, $X=\{s\},$ and ${\cal
T}=\emptyset $, the claim follows.  $\bullet$\endproof \medskip

For $s\in S$, let $V_s:=\{v\in V-s:  \ sv\not \in F_0\}$. Note that $V_s\in\mathcal{H}_0$ for $s\in S$. Let a
set-function $p_1$ on $V$ be defined as follows.

\begin{equation}p_1(V'):= \begin{cases}
 m_S(s) & \text{if $V'=V_s$ for some $s\in S$}\\
p_0(V') & \text{otherwise.}\end{cases} \label{eq:p1def}\end{equation}

The definition of $p_1$ implies that $p_1(V')$ can be positive only if $V'\in\mathcal{H}_0$.

\begin{claim}\label{claim:emeltp0} $p_1(V_s)\geq p_0(V_s)$ holds for every $s\in
S$.  \end{claim}

\proof{Proof.} Consider first the case when $V_s\cap T=\{y\}$ for some $y\in
T$.  By applying \eqref{eq:ftgs} to $X=S-s$, $Y=\{y\}$, and ${\cal
T}= \emptyset $, we get $$m_T(y) - \vert S-s\vert +d_{H_0}(y) = m_T(y)
- d_{G_0}(S-s,y)\leq m_S(s).$$ By applying \eqref{eq:ftgs} to $X=S-s$, $Y=\emptyset $, and ${\cal T}= \{y\}$, we get $p_T(y) -
r_S(S-s)\leq m_S(s)$, from which $$m_S(s) \geq \max \{ p_T(y) -
r_S(S-s), \ m_T(y) - \vert S-s\vert + d_{H_0}(y) \} = p_0(V_s).$$

Second, assume that $\vert V_s\cap T\vert \geq 2.$ By applying
\eqref{eq:ftgs} to $X=S-s$, $Y=\emptyset $, and ${\cal T}=
\{V_s\cap T\}$, we get $$p_0(V_s)= p_T(V_s\cap T) -r_S(S-s) \leq
m_S(s).\qquad \bullet$$\endproof

\begin{claim}\label{claim:ST-crossing} The set-function $p_1$ is positively
$ST$-crossing supermodular.  \end{claim}

\proof{Proof.} It follows from Claim \ref{claim:emeltp0} that $p_1$ arises from
$p_0$ by increasing its values on sets $V_s$ $(s\in S)$.  Let
$V'\subset V$ be a set which is $ST$-crossing with $V_s$ (in
particular, $V'$ and $V_s$ are not comparable). Then $S\not \subseteq V_s\cup V'$ and hence $V'\cap S\subseteq V_s\cap
S$.  Therefore $V'\cap T \not \subseteq V_s\cap T$, that is, there is
an element $t\in (V'-V_s)\cap T$.  Since $st$ is an arc of
$\overrightarrow {H_0}$ entering $V'$, we conclude that $p_1(V')=0$,
implying that $p_1$ is indeed positively $ST$-crossing supermodular.
$\bullet$\endproof \medskip

Let $\nu $ denote the maximum total $p_1$-value of $ST$-independent
sets.

\begin{lemma}$\nu =\gamma $. \end{lemma}

\proof{Proof.} Since the family ${\cal L}=\{V_s:s\in S\}$ is $ST$-independent,
$\nu \geq \widetilde p_1({\cal L}) = \widetilde m_S(S)=\gamma .$
Suppose indirectly that $\nu > \gamma $ and let $\cal I$ be an
$ST$-independent family for which $\widetilde p_1({\cal I})=\nu $. We
can assume that $\vert {\cal I}\vert $ is minimal in which case
$p_1(V')>0$ for each $V'\in \cal I$.

\begin{claim}\label{claim:no-two} There are no two $T$-intersecting members $V_1$
and $V_2$ of $\cal I$ for which $p_1(V_i)= p_0(V_i)$ $(i=1,2)$.  \end{claim}

\proof{Proof.} Suppose indirectly the existence of such $T$-intersecting
members $V_1$ and $V_2$ of $\cal I$.  Since $\cal I$ is
$ST$-independent, we must have $S\subseteq V_1\cup V_2$ and hence
$p_0(V_1\cup V_2)=0$.  Since $p_0$ is positively $T$-intersecting
supermodular,

$$p_1(V_1)+p_1(V_2) = p_0(V_1) + p_0(V_2) \leq $$ $$p_0(V_1\cap V_2) +
p_0(V_1\cup V_2) = p_0(V_1\cap V_2)\leq p_1(V_1\cap V_2).$$

Now ${\cal I}'={\cal I} -\{V_1,V_2\} + \{V_1\cap V_2\}$ is also
$ST$-independent and $\widetilde p_1({\cal I}')\geq \widetilde
p_1({\cal I})$, but we must have equality by the optimality of $\cal
I$, contradicting the minimality of $\vert {\cal I}\vert $. $\bullet$\endproof

\medskip

We say that a member $V'\in {\cal I}$ is of Type I \ if $V'=X_t+t$ for
some $t\in T$ and $X_t\subseteq S$ and
$$p_1(X_t+t)=p_0(X_t+t)=m_T(t)-\vert X_t\vert + d_{H_0}(t) >
p_T(t)-r_S(X_t).$$ Let ${\cal I}_1$ \ ($\subseteq {\cal I}$) denote
the family of sets of Type I. Claim \ref{claim:no-two} implies that if
$X_1+t_1\in {\cal I}_1$ and $X_2+t_2\in {\cal I}_1$ for which
$X_1+t_1\not =X_2+t_2$ \ ($X_i\subseteq S, \ t_i\in T$), \ then
$t_1\not =t_2$.  Let $$ \hbox{ $Y:=\{t\in T:$ there is a member
$X_t+t\in {\cal I}_1\}$.  }\ $$ Note that $\vert Y\vert =\vert {\cal
I}_1\vert $.

We say that a member $V'\in {\cal I}$ is of Type II \ if $$p_1(V') =
p_0(V') = p_T(V'\cap T) -r_S(V'\cap S).$$ Let ${\cal
I}_2=\{V_1,V_2,\cdots ,V_q\}$ ($\subseteq {\cal I}$) denote the family
of sets of Type II.  \ Let $$ \hbox{ ${\cal T}:=\{T_1,\dots ,T_q\}$
where $T_i:=V_i\cap T$ for $i=1,\dots ,q.$ }\ $$ Since $p_1(V_i)>0$,
the members of $\cal T$ are non-empty.  Furthermore, Claim
\ref{claim:no-two} implies that $\cal T$ is a subpartition of $T-Y$.

Let ${\cal I}_3:= {\cal I}- ({\cal I}_1 \cup {\cal I}_2)$.  The
members of ${\cal I}_3$ are called of Type III.  Then each member $V'$
of ${\cal I}_3$ is of form $V'=V_s$ for some $s\in S$ such that
$m_S(s)=p_1(V') > p_0(V')$.  Let $$X:=\{s\in S:  V_s\in {\cal
I}_3\}.$$

It follows from the definitions that ${\cal I}_1, {\cal I}_2$, and
${\cal I}_3$ form a partition of $\cal I$.

\begin{claim}\label{claim:XXt} Let $t\in Y$ and $X_t+t\in {\cal I}_1$.  Then
$X\subseteq X_t$.  \end{claim}

\proof{Proof.} Suppose indirectly that there is an element $s\in X-X_t$.  By
the $ST$-independence of the sets $X_t+t$ and $V_s$, the element $t$
cannot be in $V_s$.  Therefore the arc $st$ belongs to
$\overrightarrow {F_0}$.  Since $st$ enters $X_t+t$, we have
$p_1(X_t+t)=0$, a contradiction.  $\bullet$\endproof

\begin{claim}\label{claim:Xt} \ $\sum [\vert X_t\vert -d_{H_0}(t) :  t\in Y] \geq
d_{G_0}(X,Y)$.  \end{claim}

\proof{Proof.} What we prove is that $\vert X_t\vert -d_{H_0}(t) \geq
d_{G_0}(X,t)$ for $t\in Y$ and $X_t+t\in {\cal I}_1$.  Since no arc of
$\overrightarrow {H_0}$ enters $X_t+t$ and since $X\subseteq X_t$ by
Claim \ref{claim:XXt}, we have $$ \vert X_t\vert -d_{H_0}(t) = \vert
X_t\vert -d_{H_0}(X_t,t) = d_{G_0}(X_t,t) \geq d_{G_0}(X,t),$$ as
required.  $\bullet$\endproof

\begin{claim}\label{claim:Xresze} $X\cup\Gamma_{H_0}(T_i)\subseteq V_i\cap S$ holds for each $i=1,\dots,q$.\end{claim}

\proof{Proof.} As $V_i\in\mathcal{H}_0$, we have $\Gamma_{H_0}(T_i)\subseteq V_i\cap S$. If, indirectly, there is an $s\in X-V_i$, then the
$ST$-independence of $V_s$ and $V_i$ implies that $V_s\cap V_i\cap
T=\emptyset $. In this case, an element $t\in V_i\cap T$ cannot be in
$V_s$ implying that $st\in \overrightarrow {F_0}$.  But then $p_1(V_i)
= 0$, contradicting the property $p_1(V')>0$ for each $V'\in {\cal
I}.$ $\bullet$\endproof

\medskip

Recall that ${\cal T}$ is a subpartition of $T - Y$.  This and the
last two claims imply $$ \gamma < \nu = \widetilde p_1({\cal I}) =
\widetilde p_1({\cal I}_1) + \widetilde p_1({\cal I}_2) + \widetilde
p_1({\cal I}_3)=$$ $$\sum [m_T(t) - \vert X_t\vert + d_{H_0}(t) :
X_t+t\in {\cal I}_1] + \sum _{i=1}\sp q [p_T(T_i) - r_S(V_i\cap S)] +
\sum [m_S(s):  V_s\in {\cal I}_3] \leq $$ $$ \sum [m_T(t):  X_t+t\in
{\cal I}_1] - d_{G_0}(X,Y) + \sum _{i=1}\sp q [p_T(T_i) - r_S(X\cup
\Gamma _{H_0}(T_i))] + \widetilde m_S(X)=$$ $$\widetilde m_T(Y) -
d_{G_0}(X,Y) + \sum _{i=1}\sp q [p_T(T_i) - r_S(X\cup \Gamma
_{H_0}(T_i))] + \widetilde m_S(X),$$ in a contradiction with
\eqref{eq:ftgs}, completing the proof of the lemma.  $\bullet $ $\bullet
$\endproof \medskip

\begin{claim} \label{claim:canuse}
If $p_1(V')$ is positive, then $\overrightarrow {G_0}$ covers $V'$.
\end{claim}
\proof{Proof.}
As already observed after \eqref{eq:p1def}, $V'\in\mathcal{H}_0$. Assume to the contrary that $\overrightarrow {G_0}$ does not cover $V'$. As $G_0$ denotes the bipartite complement of $H_0$, this is only possible if $V'\cap T=\emptyset$ or $S\subseteq V'$.

If $V'=V_s$ for some $s\in S$, then $s\notin V'$, hence $V'\cap T=\emptyset$. This means that $st\in F_0$ for each $t\in T$. By applying \eqref{eq:ftgs} to $X=\{s\}$, $Y=T$ and $\mathcal{T}=\emptyset$, we get $p_1(V')=m_S(s)\leq 0$, a contradiction.

Therefore, we must have $p_1(V')=p_0(V')$. As $p_0$ was defined to be $0$ for sets not intersecting $T$, we can assume that $S\subseteq V'$ holds. If $p_0(V')=p_T(V'\cap T)-r_S(V'\cap S)$, then \eqref{eq:ftgs}, when applied to $Y=\emptyset$, $X=S$ and $\mathcal{T}=\{V'\cap T\}$, gives $p_0(V')=p_T(V'\cap T)-r_S(S)\leq 0$, a contradiction. Therefore $p_0(V')$ is defined by the first line of \eqref{eq:p0def}. Hence $V'=S+y$ for some $y\in T$ and $p_0(V')=m_T(y)-|S|+d_{H_0}(y)$. Now \eqref{eq:ftgs}, when applied to $Y=\{y\}$, $X=S$ and $\mathcal{T}=\emptyset$, gives $p_0(V')=m_T(y)-|S|+d_{H_0}(y)=m_T(y)-d_{G_0}(S,y)\leq 0$, thus leading to a contradiction again.
$\bullet$\endproof

\medskip

By Claim~\ref{claim:canuse}, Theorem \ref{thm:Frank-Jordan} can be applied to $p_1$. This means that there is a digraph $D=(V,A)$ on $V$
with $\nu = \gamma $ $ST$-arcs that covers $p_1$, that is, $\varrho
_D(V')\geq p_1(V')$ for every subset $V'\subseteq V$.  Let $G=(S,T;E)$
denote the undirected bipartite graph underlying $D$.

\begin{claim}$d_G(s)=m_S(s)$ for every $s\in S$ and $d_G(t)=m_T(t)$ for
every $t\in T$.\end{claim}

\proof{Proof.} Since $d_G(s)=\delta _D(s)\geq \varrho _D(V_s)\geq
p_1(V_s)=m_S(s)$ for every $s\in S$, we have \ $\widetilde
m_S(S)=\vert E\vert = \sum [ d_G(s):s\in S] \geq \widetilde m_S(S)$,
from which $d_G(s)=m_S(s)$ follows for every $s\in S$.

Since $d_G(t)=\varrho _D(t) \geq \varrho _D(\Gamma _{H_0}(t)+t) \geq
p_0(\Gamma _{H_0}(t) +t)\geq m_T(t)$ for every $t\in T$, we have
$\widetilde m_T(T)=\vert E\vert = \sum [ d_G(t):t\in T] \geq
\widetilde m_T(T)$, from which $d_G(t)=m_T(t)$
follows for every $t\in T$.  $\bullet$\endproof

\begin{claim}\label{claim:G+simple} The bigraph $G\sp +=(S,T;E+F_0)$ is simple.
\end{claim}

\proof{Proof.} The minimality of $D$ implies that each arc of $D$ enters a
subset $V'$ with $p_1(V')>0$.  Since $p_1(V')$ can be positive only if
no arc of $\overrightarrow {H_0}$ enters $V'$, we can conclude that no
edge of $G$ is parallel to an edge of $H_0$.

Suppose indirectly that there are two parallel edges $e$ and $e'$ of
$G$ connecting $s$ and $t$ for some $s\in S$ and $t\in T$.  Let
$X:=\Gamma_{H_0}(t)$.  Then $p_1(X+t)\geq m_T(t) =\varrho
_D(t)$.  For $V'=X+s+t$, we have $\varrho _D(t) -2 \geq \varrho
_D(V')\geq p_1(V')\geq p_1(X+t) -1 \geq m_T(t)-1=\varrho _D(t)-1$, a
contradiction.  $\bullet$\endproof

\begin{claim}\label{claim:Mfed} $r_S(\Gamma _{G\sp +}(Y))\geq p_T(Y)$ for every
subset $Y\subseteq T$.  \end{claim}

\proof{Proof.} Let $X:= \Gamma _{G\sp +}(Y)$ and $V':=Y\cup X$.  Then $0=
\varrho _D(V')\geq p_1(V') \geq p_T(Y) -r_S(X)$, from which the claim
follows.  $\bullet$\endproof \medskip

We conclude that $G$ meets all the requirements of the theorem, and
the proof is complete.  $\bullet $ $\bullet $ $\bullet $ \endproof \medskip

\subsection{Variations}

\subsubsection{Degree-specification only on \texorpdfstring{$S$}{S}}

With the proof technique used above, one can derive the following
variation where the degrees are specified only for the nodes in $S$.
Namely, the definition of $p_0$ in \eqref{eq:p0def} should be modified
as follows.

\begin{equation}p_0(V'):= \begin{cases}
 p_T(Y)-r_S(X) & \text{if $V'=X\cup Y\in {\cal H}_0$,
\ $X\subseteq S$, \ $Y\subseteq T$}\\ 0 & \text{otherwise.} \end{cases}
\label{eq:p0defb} \end{equation}

\begin{theorem}\label{thm:ms.nov} We are given a simple bigraph $H_0=(S,T;F_0)$,
a matroid $M_S=(S,r_S)$, a positively intersecting supermodular
function $p_T$ on $T$, and a degree-specification $m_S$ on $S$ for
which $\widetilde m_S(S)=\gamma $. There is a bigraph $G=(S,T;E)$
fitting $m_S$ for which $G\sp +=G+H_0$ is simple and $M_S$-covers
$p_T$ if and only if \begin{equation}m_S(s) +d_{H_0}(s)\leq \vert T\vert \
\hbox{for every $s\in S$}\ \label{eq:ujterm.msb} \end{equation} and \begin{eqnarray}&\widetilde m_S(X) + \sum _{i=1}\sp q [p_T(T_i) -
r_S(X\cup \Gamma _{H_0}(T_i)) ] \leq \gamma& \nonumber \\ &\hbox{whenever
$X\subseteq S$ and \ ${\cal T}=\{T_1,\dots ,T_q\}$ a subpartition of
$T$. }\ & \label{eq:ms.nov} \end{eqnarray} \end{theorem}

One reason why we do not go into the details is that the proof is
quite similar to (and, in fact, slightly simpler than) the proof of
Theorem \ref{thm:msmt.nov}.  Another reason is that, in a forthcoming work
\cite{Berczi-Frank17}, we solve a common generalization of Theorems
\ref{thm:msmt.nov} and \ref{thm:ms.nov} where, instead of
degree-specifications, there are both upper and lower bounds for the
degrees of all nodes in $S\cup T$.

\subsubsection{Fully supermodular \texorpdfstring{$p_T$}{pT}}

In the special case when $p_T\equiv 0$, it suffices to require the
inequality in \eqref{eq:ftgs} only for the empty $\cal T$, in which case
Theorem \ref{thm:msmt.nov} reduces to the following classic result (which
actually holds for non-simple bigraphs, too).

\begin{theorem}[Ore \cite{Ore56}] \label{thm:Ore.G0} A simple bigraph
$G_0=(S,T;E_0)$ has a subgraph fitting a degree-specification
$(m_S,m_T)$ with $\widetilde m_S(S)=\widetilde m_T(T)=\gamma $ if and
only if \begin{equation}\widetilde m_S(X) + \widetilde m_T(Y) -d_{G_0}(X,Y) \leq
\gamma \ \hbox{ whenever $X\subseteq S, \ Y\subseteq T$.}\
\label{eq:Ore.felt} \end{equation} \end{theorem}

The content of the next result is that the condition in Theorem
\ref{thm:msmt.nov} can also be simplified when $p_T$ is fully
supermodular.

\begin{theorem}\label{thm:msmt.fully} We are given a simple bigraph
$H_0=(S,T;F_0)$, a matroid $M_S=(S,r_S)$, a fully supermodular
function $p_T$ on $T$, and a degree-specification $m_V=(m_S,m_T)$ for
which $\widetilde m_S(S)=\widetilde m_T(T)=\gamma $. There is a
bigraph $G=(S,T;E)$ fitting $m_V$ for which $G\sp +=G+H_0$ is simple
and $M_S$-covers $p_T$ if and only if \eqref{eq:Ore.felt} holds and
\begin{eqnarray}&\widetilde m_S(X)+ \widetilde m_T(Y) -
d_{G_0}(X,Y) + p_T(T_0) - r_S(X\cup \Gamma _{H_0}(T_0)) \leq \gamma & \nonumber\\
&\hbox{whenever $Y\subseteq T$, $X\subseteq S$, $T_0\subseteq T-Y$, }\
& \label{eq:ftgs.fully} \end{eqnarray} where $G_0$ is the
bipartite complement of $H_0$.  \end{theorem}

\proof{Proof.} Conditions \eqref{eq:Ore.felt} and \eqref{eq:ftgs.fully} correspond
to the special cases of Condition \eqref{eq:ftgs} when $\vert {\cal
T}\vert =0$ and $\vert {\cal T}\vert =1$, respectively.  Therefore
their necessity was proved earlier.  To see sufficiency, by Theorem
\ref{thm:msmt.nov} it suffices to show that \eqref{eq:ftgs} holds in
general.  Suppose, indirectly, that there are $X,$ $Y$, and ${\cal T}$
violating \eqref{eq:ftgs}.  Assume that $\vert {\cal T}\vert $ is
minimal.  Then \eqref{eq:Ore.felt} and \eqref{eq:ftgs.fully} imply that
$\vert {\cal T}\vert \geq 2$.  Let $T_1,T_2$ be two members of $\cal
T$.  Since $$p_T(T_1\cup T_2) -r_S(X\cup \Gamma _{H_0}(T_1\cup T_2))
\geq p_T(T_1) + p_T(T_2) - r_S(X\cup \Gamma _{H_0}(T_1)) - r_S(X\cup
\Gamma _{H_0}(T_2)) ,$$ the unchanged sets $X,Y$ and the partition
${\cal T}'$ obtained from $\cal T$ by replacing $T_1$ and $T_2$ with
the single set $T_1\cup T_2$ also violate \eqref{eq:ftgs}, contradicting
the minimal choice of $\cal T$.  $\bullet$\endproof

\medskip

It is worth formulating Theorem \ref{thm:msmt.fully} in the special case
when $H_0$ has no edges.

\begin{corollary}\label{cor:csak.matroid} We are given a matroid $M_S=(S,r_S)$,
a fully supermodular function $p_T$ on $T$, and a degree-specification
$m_V=(m_S,m_T)$ for which $\widetilde m_S(S)=\widetilde m_T(T)=\gamma
$. There is a simple bigraph $G=(S,T;E)$ fitting $m_V$ and
$M_S$-covering $p_T$ if and only if \begin{equation}\widetilde m_S(X) + \widetilde
m_T(Y) -\vert X\vert \vert Y\vert \leq \gamma \ \hbox{ whenever
$X\subseteq S, \ Y\subseteq T$}\ \label{eq:Ore.felt.0} \end{equation} and \begin{equation}
\widetilde m_S(X) + \widetilde m_T(Y) - \vert X\vert \vert Y\vert +
p_T(T_0) - r_S(X) \leq \gamma \ \hbox{whenever $Y\subseteq T$,
$X\subseteq S$, $T_0\subseteq T-Y$.  \ }\ \label{eq:csak.matroid} \end{equation}
If, in addition, $p_T$ is monotone non-decreasing, then $T_0$ in
\eqref{eq:csak.matroid} can be chosen to be $T_0=T - Y$, that is, \begin{equation}
\widetilde m_S(X)+ \widetilde m_T(Y) - \vert X\vert \vert Y\vert +
p_T(T-Y) - r_S(X) \leq \gamma \ \hbox{whenever $X\subseteq S$,
$Y\subseteq T.$  }\ \label{eq:csak.matroid.mon} \end{equation} \end{corollary}

\proof{Proof.} The first part is a special case of Theorem \ref{thm:msmt.fully}.
When $p_T$, in addition, is monotone non-decreasing in the second
part, we can choose $T_0$ in \eqref{eq:csak.matroid} as large as
possible, that is, $T_0 = T - Y$.  $\bullet$\endproof

\medskip 

\section{Matroidal max term rank}

Let ${\cal G}(m_S,m_T)$ denote the set of simple bigraphs $G=(S,T;E)$
fitting a degree-specification $(m_S,m_T)$ with $\widetilde
m_S(S)=\widetilde m_T(T)=\gamma $. It follows from Theorem
\ref{thm:Ore.G0} that ${\cal G}(m_S,m_T)$ is non-empty if and only if
\eqref{eq:Ore.felt.0} holds.  In \cite{Berczi-Frank16a} (Theorem 26),
we reformulated Ryser's classic max term rank formula in graph
theoretic language.

\begin{theorem}[Ryser] \label{thm:Ryser} Let $\ell\leq \vert T\vert $ be an
integer.  Suppose that ${\cal G}(m_S,m_T)$ is non-empty.  Then ${\cal
G}(m_S,m_T)$ has a member $G$ with matching number $\nu (G)\geq \ell$
if and only if \begin{equation}\widetilde m_S(X) + \widetilde m_T(Y) -\vert X\vert
\vert Y\vert + (\ell - \vert X\vert -\vert Y\vert ) \leq \gamma \
\hbox{ whenever}\ X\subseteq S, \ Y\subseteq T. \label{eq:Ryser} \end{equation}
Moreover, \eqref{eq:Ryser} holds if the inequality in it is required
only when $X$ consists of the $i$ largest values of $m_S$ and $Y$
consists of the $j$ largest values of $m_T$ $(i=0,1,\dots ,\vert
S\vert , \ j=0,1,\dots ,\vert T\vert )$.  \end{theorem}

We keep using graph terminology, but the original expression (max term
rank) of Ryser is retained.  Our present goal is to extend Ryser's
theorem in two directions.  In the augmentation version an initial
bigraph is to be augmented while in the matroidal form the matching is
expected to cover a basis of a matroid $M_S$ on $S$ and a basis of
matroid $M_T$ on $T$.  Actually, we shall integrate the two
generalizations into one single framework.

In what follows, $M_S=(S,r_S)$ and $M_T=(T,r_T)$ will be matroids of
rank $\ell$.  In \cite{Berczi-Frank16a}, the complementary
set-function $p$ of a set-function $b$ was defined by
$p(Y):=b(S)-b(S-Y)$.  Clearly, $b$ is submodular if and only if $p$ is
supermodular, and $p$ is monotone non-decreasing if and only if $b$ is
monotone non-decreasing.  The complementary function $p_T$ of the rank
function $r_T$ of $M_T$ is called the {\bf co-rank function} of $M_T$.
It can easily be shown that $p_T(Y)= \min \{\vert Y\cap B\vert :  B$ a
basis of $M_T\}$.

The following extension of Edmonds' matroid intersection theorem
\cite{Edmonds70} will be used.  For notational convenience, the
bipartite graph in the theorem is denoted by $G\sp +$.

\begin{theorem}[Brualdi, \cite{Brualdi70b}] \label{thm:Brualdi} Let $G\sp
+=(S,T;E\sp +)$ be a bigraph with a matroid $M_S=(S,r_S)$ on $S$ and
with a matroid $M_T=(T,r_T)$ on $T$ for which $r_S(S)=r_T(T)=\ell.$
There is a matching of $G\sp +$ covering bases of $M_S$ and $M_T$ if
and only if \begin{eqnarray}&r_S(X') + r_T(Y') \geq \ell& \nonumber \\ & \hbox{whenever $X'\cup Y'$ hits every edge of $G\sp +$ \
$(X'\subseteq S, \ Y'\subseteq T)$.}& \label{eq:termrank.Bru} \end{eqnarray} \end{theorem}

We need the following equivalent version of Theorem \ref{thm:Brualdi}.

\begin{lemma}\label{lem:Brualdib} Let $G\sp +=(S,T;E\sp +)$ be a bigraph.  Let
$M_S$ be a matroid on $S$ with rank function $r_S$ and $M_T$ a matroid
on $T$ with co-rank function $p_T$ for which $r_S(S)=p_T(T)=\ell$.
There is a matching of $G\sp +$ covering bases of $M_S$ and $M_T$ if
and only if \begin{equation}r_S(\Gamma _{G\sp +}(Y)) \geq p_T(Y) \ \hbox{for
every}\ Y\subseteq T. \label{eq:termrank.Brub} \end{equation} \ \end{lemma}

\proof{Proof.} The necessity is straightforward.  The sufficiency follows from
Theorem \ref{thm:Brualdi} once we show that \eqref{eq:termrank.Bru} holds.
Since $X'\cup Y'$ hits every edge of $G\sp +$, for $Y:=T-Y'$ we have
$\Gamma _{G\sp +}(Y) \subseteq X'$.  Therefore \eqref{eq:termrank.Brub}
implies that $r_S(X')\geq r_S(\Gamma _{G\sp +}(Y))\geq p_T(Y)= r_T(T)
- r_T(Y') = \ell - r_T(Y')$ and hence \eqref{eq:termrank.Bru} indeed
holds.

\begin{theorem}\label{thm:Ryser.gen} We are given a simple bigraph
$H_0=(S,T;F_0)$, a matroid $M_S=(S,r_S)$ and a matroid $M_T=(T,r_T)$
with $r_S(S)=r_T(T)=\ell,$ and a degree-specification $m_V=(m_S,m_T)$ on $V:=S\cup T$ for which $\widetilde m_S(S)=\widetilde m_T(T)=\gamma $. There is a
bigraph $G=(S,T;E)$ fitting $m_V$ for which $G\sp +=G+H_0$ is simple
and includes a matching covering a basis of $M_S$ and a basis of $M_T$
if and only if \eqref{eq:Ore.felt} holds and \begin{eqnarray}& \widetilde m_S(X) + \widetilde m_T(Y) - d_{G_0}(X,Y)
+ \ell - r_S(X') -r_T(Y') \leq \gamma & \nonumber \\ & \text{whenever $X\subseteq
X'\subseteq S$, $Y\subseteq Y'\subseteq T$, and $X'\cup Y'$ hits all the
edges of $H_0$, }\ & \label{eq:Ryser.gen} \end{eqnarray}
where $G_0$ is the bipartite complement of $H_0$.  \end{theorem}

\proof{Proof.} Necessity.  Suppose that the requested bigraph $G$ and its
$\ell$-element matching $M$ exist.  The number of edges of $G$ with at
least one end-node in $X\cup Y$ is at least $\widetilde m_S(X) +
\widetilde m_T(Y) - d_{G_0}(X,Y)$.  The number of edges in $M$ with at
least one end-node in $X'\cup Y'$ is at most $r_S(X') + r_T(Y')$.
Therefore $M$ has at least $\ell - r_S(X') - r_T(Y')$ elements
connecting $S-X'$ and $T-Y'$.  But these elements must be in $E$ since
$X'\cup Y'$ hits all edges of $H_0$.  Therefore the total number of
edges of $G$ is at least $\widetilde m_S(X) + \widetilde m_T(Y) -
d_{G_0}(X,Y) + \ell - r_S(X') - r_T(Y')$, and \eqref{eq:Ryser.gen}
follows.

Sufficiency.  Let $p_T$ denote the co-rank function of $M_T$, that is,
$p_T(Z)=\ell -r_T(T-Z)$ for $Z\subseteq T$. Note that $p_T$ is fully supermodular. 

\begin{claim}\label{claim:origtrue}
Condition \eqref{eq:ftgs.fully} is satisfied.
\end{claim}
\proof{Proof.}
For the present $p_T$, Condition \eqref{eq:ftgs.fully} requires
\begin{equation}
\widetilde m_S(X)+ \widetilde m_T(Y) - d_{G_0}(X,Y) + \ell - r_T(T-T_0) - r_S(X\cup\Gamma_{H_0}(T_0)) \leq \gamma. \label{eq:Ryser.gen.regi}
\end{equation}
No sets $X\subseteq S$, $Y\subseteq T$, and $T_0\subseteq T-Y$ can violate this inequality since then, by letting $Y':= T-T_0$ and $X':=X\cup\Gamma_{H_0}(T_0)$, the quadruple $(X,\ Y,\ X',\ Y')$ would violate \eqref{eq:Ryser.gen}. 
$\bullet$\endproof

\medskip

By Theorem \ref{thm:msmt.fully}, there is a bigraph $G$
fitting $m_V$ for which $G\sp +=G+H_0$ is simple and $M_S$-covers
$p_T$.  The latter property, by definition, means that
\eqref{eq:termrank.Brub} holds, and therefore Lemma \ref{lem:Brualdib}
implies that $G\sp +$ has a requested matching.$\bullet$ $\bullet$\endproof 

\medskip

When $m_V\equiv 0$ and $\gamma =0$, it suffices to require
\eqref{eq:Ryser.gen} only for $X=Y=\emptyset $ in which case it
transforms to \begin{equation} \ell - r_S(X') -r_T(Y')
\leq 0\ \text{whenever $X'\subseteq S$, $Y'\subseteq T$, and $X'\cup Y'$
hits all the edges of $H_0$, }\label{eq:Ryser.genx} \end{equation} which is the same as \eqref{eq:termrank.Bru}.
In other words, Theorem \ref{thm:Ryser.gen} may be considered as a
straight generalization of Brualdi's theorem.

The content of the next corollary is that in the special case of
Theorem 8 when $F_0=\emptyset $ it suffices to require \eqref{eq:Ryser.gen} only in a simplified form.

\begin{corollary}\label{cor:Ryser.matroid} Let $S$ and $T$ be two disjoint sets
and $(m_S,m_T)$ a degree-specification on $S\cup T$ for which
$\widetilde m_S(S)=\widetilde m_T(T)=\gamma $ and ${\cal G}(m_S,m_T)$ is non-empty, that is,
\eqref{eq:Ore.felt.0} holds.  Let $M_S=(S,r_S)$ and $M_T=(T,r_T)$ be
matroids for which $r_S(S)=r_T(T)=\ell.$ There is a simple bigraph
$G=(S,T;E)$ fitting $(m_S,m_T)$ that includes a matching covering bases of
$M_S$ and $M_T$ if and only if \begin{equation}\widetilde m_S(X) + \widetilde
m_T(Y) - \vert X\vert \vert Y\vert + \ell - r_S(X) - r_T(Y) \leq
\gamma \label{eq:Ryser.matroid} \end{equation} holds for every $X\subseteq S$ and
$ Y\subseteq T$.\end{corollary}
\proof{Proof.}
Consider Theorem \ref{thm:Ryser.gen} in the special case when $F_0=\emptyset $.
Then the bipartite complement $G_0$ of $H_0$ is a complete bigraph and hence $d_{G_0}(X,Y)= \vert X\vert
\vert Y\vert $. Therefore Condition \eqref{eq:Ore.felt.0} requested in the corollary is the same as Condition \eqref{eq:Ore.felt} requested in Theorem \ref{thm:Ryser.gen}. Furthermore \eqref{eq:Ryser.matroid} is
the special case of \eqref{eq:Ryser.gen} when $X'=X$ and $Y'=Y$, and hence \eqref{eq:Ryser.matroid} is
necessary.

We claim, conversely, that \eqref{eq:Ryser.matroid} implies \eqref{eq:Ryser.gen}.  Indeed, if $X\subseteq X'$ and $Y\subseteq Y'$ violate the inequality in \eqref{eq:Ryser.gen}, then the monotonicity of matroid rank functions imply that $\gamma < \widetilde
m_S(X) + \widetilde m_T(Y) - \vert X\vert \vert Y\vert + \ell -
r_S(X') - r_T(Y') \leq \widetilde m_S(X) + \widetilde m_T(Y) - \vert
X\vert \vert Y\vert + \ell - r_S(X) - r_T(Y)$, contradicting \eqref{eq:Ryser.matroid}.
Therefore the requested bigraph exists by Theorem \ref{thm:Ryser.gen}.
$\bullet$\endproof

\medskip

Note that the inequalities in Conditions \eqref{eq:Ore.felt.0} and
\eqref{eq:Ryser.matroid} can be integrated into the following single
form:

\begin{equation}\widetilde m_S(X) + \widetilde m_T(Y) - \vert X\vert \vert Y\vert
+ (\ell - r_S(X) - r_T(Y))\sp + \leq \gamma . \end{equation}

By specializing Theorem \ref{thm:Ryser.gen} to the case when $M_S$ and
$M_T$ are $\ell$-uniform matroids on $S$ and $T$, respectively, one
obtains the following.

\begin{corollary}\label{cor:Ryser.novel} We are given a simple bigraph
$H_0=(S,T;F_0)$, an integer $\ell$, and a degree-specification
$(m_S,m_T)$ for which $\widetilde m_S(S)=\widetilde m_T(T)=\gamma
$. There is a bigraph $G=(S,T;E)$ fitting $(m_S,m_T)$ for which $G\sp
+=G+H_0$ is simple and includes an $\ell$-element matching if and only
if \eqref{eq:Ore.felt} holds and \begin{eqnarray}&\widetilde m_S(X)+ \widetilde
m_T(Y) - d_{G_0}(X,Y) + \ell - \vert X'\vert
-\vert Y'\vert \leq \gamma& \nonumber \\ &\text{whenever $X\subseteq
X'\subseteq S$, $Y\subseteq Y'\subseteq T$, and $X'\cup Y'$ hits all the
edges of $H_0$, } & \label{eq:Ryser.novel} \end{eqnarray} where $G_0$ is
the bipartite complement of $H_0$.\end{corollary}
\proof{Proof.}
Consider Theorem \ref{thm:Ryser.gen} in the special case when $M_S$ and $M_T$ are
$\ell$-uniform matroids on $S$ and on $T$, respectively.  Since
matroid rank functions are subcardinal, \eqref{eq:Ryser.novel} is implied by \eqref{eq:Ryser.gen} and hence \eqref{eq:Ryser.novel} is necessary.

We claim, conversely, that \eqref{eq:Ryser.novel} implies \eqref{eq:Ryser.gen}, that is, $\alpha + \ell
- r_S(X') - r_T(Y') \leq \gamma $ where $\alpha := \widetilde m_S(X) +
\widetilde m_T(Y) - d_{G_0}(X,Y)$.  Indeed, if $\max\{\vert X'\vert
,\vert Y'\vert \} \leq \ell$, then $ \alpha + \ell - r_S(X') - r_T(Y')
= \alpha + \ell - \min\{ \ell, \vert X'\vert \} - \min\{\ell, \vert
Y'\vert \} = \alpha + \ell - \vert X'\vert - \vert Y'\vert \leq \gamma
$, where the last inequality follows by \eqref{eq:Ryser.novel}.  If $\max\{\vert X'\vert
,\vert Y'\vert \} > \ell$, and, say, $\vert X'\vert >\ell$, then
$\alpha + \ell - r_S(X') - r_T(Y') = \alpha +\ell - \min\{ \ell, \vert
X'\vert \} - \min\{\ell, \vert Y'\vert \} \leq \alpha + \ell - \ell -
0 = \alpha \leq \gamma $, where the last inequality follows by \eqref{eq:Ore.felt}.
Therefore the requested bigraph exists by Theorem \ref{thm:Ryser.gen}.
$\bullet$\endproof

\medskip

\noindent {\bf Acknowledgements} \ We are grateful to Richard Brualdi
for the valuable information on the topic.  Special thanks are due to
Zolt\'an Szigeti for carefully checking the details of a first draft.
The two anonymous referees provided several useful comments.  We
gratefully acknowledge their great contribution.

The research was supported by the Hungarian Scientific Research Fund -
OTKA, No K109240.  The work of the first author was financed by a
postdoctoral fellowship provided by the Hungarian Academy of Sciences.

\medskip 

\end{document}